\documentclass[11pt]{article}\usepackage[english]{babel}
\usepackage{amsfonts,amssymb,amsmath,amsthm,latexsym}
\usepackage[usenames]{color}
\usepackage{inputenc}
\usepackage{enumerate}
\newtheorem{teor}{Theorem}[section]

\newtheorem{defi}[teor]{Definition}

\newtheorem{ejm}[teor]{Example}
\newtheorem{cor}[teor]{Corollary}

\newtheorem{remark}[teor]{Remark}

\newcommand{\PC}{\mathbb{P}}

\newcommand{\Z}{\mathbb{Z}}

\newcommand{\C}{\mathbb{C}}

\newcommand{\noi}{\noindent}
\newcommand{\prova}{\noi{\bf Proof}}
\newcommand{\fin}{\hfill$\square$\\}

\newcommand{\inter}{\cap}

\DeclareMathOperator{\I}{I}

\title{On the factorization of the polar of a  plane branch}

\author{A. Hefez, M. E. Hernandes and M. F. H. Iglesias \thanks{The first two authors were partially supported by a CNPq grant, while the third author was supported by a fellowship from CAPES/Funda\c c\~ao Arauc\'aria.}}

\date{ \ }
\begin{document}

\maketitle

\begin{abstract}
In this paper we present the most complete description as possible of the factorization of the general polar of the general member of an equisingularity class of irreducible germs of complex plane curves. Our result will refine the rough description of the factorization given by M. Merle in \cite{M} and it is based on the result given by E. Casas-Alvero in \cite{C2} that describes the cluster of the singularities of such polars. By using our analysis, it will be possible to characterize all equisingularity classes 
of irreducible plane germs with $r$ characteristic exponents having the exceptional behavior that the general polar of a general curve in this equisingularity class has only irreducible components with less than $r$ characteristic exponents, generalizing a result obtained for $r=2$ in \cite{HHI1}.
\end{abstract}

\section{Introduction}
The study of the polar of a germ of plane curve is a classical subject and the topological or equisingular classification of polars of equisingular complex plane curve germs is still an open problem. These objects have been used in classical algebraic geometry for enumeration purposes, such as Pl\"ucker formulas, and were resuscitated during the 70's in the work of B. Teissier \cite{T} for the study of families of singular hypersurfaces, being still actively studied nowadays.

Initially, it was thought that topologically equivalent germs of plane curves had topologically equivalent polar curves, which is false as shown with a simple example in \cite{P}. The topological type of the general polar of the germ of a plane curve is actually an analytic invariant of the germ. However, there are some particular invariants attached to the polars of topologically equivalent plane curve germs, namely the polar quotients, associated to the decomposition of the polar, roughly described by a Theorem of Merle in \cite{M}, as we will explicit later. Actually, in this paper we will describe topologically the complete decomposition of the general polar of the general complex plane curve germ topologically equivalent to a given irreducible one. This will be done by using the description given by E. Casas-Alvero in \cite{C2} of the cluster of the general polar of the general curve belonging to an equisingularity class of irreducible germs of complex plane curves. We explore this cluster in order to recover the characteristic exponents of each irreducible component of the general polar and the intersection multiplicities of all pairs of such components. This determines in Zariski's way the equisingularity class of the general polar of the general member of the equisingularity class of the curve. Our analysis will allow us, as a byproduct, to characterize the equisingularity classes of irreducible plane germs such that their general members have general polar that admit only irreducible components with at most one less characteristic exponent, generalizing a result obtained in \cite{HHI1} in the case of curves with two characteristic exponents.

\section{Classical results} 
A germ of an analytic plane curve at the origin of $\C^2$ is a germ of set $C=C_f=\{(x,y)\in(\mathbb{C}^2,0);\ f(x,y)=0\}$, where $f\in \mathbb{C}\{x,y\}$ is a convergent complex power series in two variables at the origin. Two such germs will be considered analytically equivalent if there is a germ of analytic diffeomorphism $\varphi$ of $(\mathbb{C}^2,0)$, also called an analytic change of coordinates, such that $\varphi(C_f)=C_g$. When the above $\varphi$ is just a homeomorphism, we say that $C_f$ and $C_g$ are topologically equivalent, or equisingular, writing, in this case, $C_f\equiv C_g$. 

From now on, we will assume that $f$ is an irreducible power series and call its associated curve $C=C_f$ a branch. After an analytic change of coordinates, if necessary, we may assume that $\I(f,x)=n<m=\I(f,y)$ and $n\nmid m$, where $\I(f,g)$ stands for the intersection multiplicity at the origin of the plane curve germs $C_f$ and $C_g$. The integer $n$ is called the multiplicity of $C$.  With such coordinates suitably chosen, it is well known that a branch $C$ admits a Newton-Puiseux parametrization of the form $(t^2,t^m)$, if $n=2$, or $(t^{n},\sum_{m\leq i <c} a_it^i)$, if $n>2$, where $c$ is some positive integer, called the conductor of $C$. Conversely, given such a parametrization, attached to it there is a well defined branch. It is also classically known that the topological, or equisingularity class of $C$ is completely determined by $n$ and the characteristic exponents $m_1,\ldots,m_r$, defined by  
$$ m_i=\min\{j;\ a_j\neq 0\ \mbox{and}\ e_{i-1}\nmid j \},$$
where $e_0=n$ and, for $k>0$, $e_k=\gcd(n,m_1,\ldots ,m_k)$ and $e_r=1$. The integer $r$ is what we call the genus of $C$. We also define the integers $d_0=1$ and $ d_i=\frac{e_{i-1}}{e_i}$, for $i=1,\ldots,r$. 

When a germ of curve is not irreducible, but reduced, Zariski has shown that its equisingularity type is determined by the equisingularity type of its branches and by their mutual intersection multiplicities.

In what follows, we will consider the set $K(n,m_1,\ldots,m_r)$ that par\-ame\-tri\-zes all Newton-Puiseux finite expansions as above with multiplicity $n$ and characteristic exponents $m_1,\ldots, m_r$.

Let $f$ be a reduced power series. The germ of curve defined by
$P_{(a:b)}(f)\colon=af_x+bf_y=0$ is the polar curve of $f$ in the
direction $(a:b)\in \PC^1$. When $(a:b)$ is a general point of $\PC^1$, we say that the associated polar $P_{(a:b)}(f)=0$ is general and we denote it shortly by $P(f)$.

In this paper we only consider the general polar of $f$ and we refer to
it simply as the polar curve.

In general, the polar curve depends upon the equation $f$ of $C_f$, however
its topological type depends only upon the analytic type of $C_f$ (see 
\cite[Theorem 7.2.10]{C3}).

The next result due to M. Merle provides a rough decomposition of $P(f)$ in 
packages of curves, not necessarily irreducible, that gives partial information about the topology of $P(f)$.

\begin{teor}[Merle \cite{M}]\label{Teorema de Merle}
Let $C_f$ be a germ of an irreducible curve with multiplicity $n$ and characteristic exponents $m_1,\ldots, m_r$. Then the general polar $P(f)$ has a decomposition of the form
$$P(f)=\gamma_1\gamma_2\cdots \gamma_r,$$ where each $\gamma_i$, not necessaritily irreducible, satisfies the following conditions:
\begin{enumerate}[{\rm i)}]
\item The multiplicity of $\gamma_i$ is given by 
$m(\gamma_i)=d_0d_1d_2\cdots d_{i-1}(d_i-1)$;
\item Each irreducible factor $\gamma_{i,j}$ of $\gamma_i$ satisfies
$$\frac{\I(\gamma_{i,j},f)}{m(\gamma_{i,j})}=\frac{1}{n} \sum_{k=1}^{i-1}(e_{k-1}-e_k) m_k+m_i.$$
\end{enumerate}
\end{teor}

Let us make some few remarks. Merle's Theorem does not describe completely the
topology of $P(f)$, because it does not describe the branches inside each package $\gamma_i$. Such branches depend upon the analytic type of $f$ and not only upon its topological type. It also does not describe the intersection multiplicities among the branches of the polar. The terms in the second conclusion are the so called polar quotients and the equality says that the branches $\gamma_{i,j}$ have contact order with $C_f$ equal to $m_i$, which implies that they have genus at least $i-1$, but they may have greater genus. 

On the other hand, Casas in \cite{C2}, determines the
equisingularity class of $P(f)$, for an $f$ corresponding to a general member of $K(n,m_1,\ldots,m_r)$ in terms of a certain weighted cluster obtained from the Enriques diagram attached to the resolution of $C_f$. 

If $r=1$, Casas in \cite{C1} describes more explicitly the factorization of $P(f)$ as follows:

Let $n$ and $m$ be two coprime natural numbers. Consider the euclidean GCD algorithm applied to the pair $n,m$:

\begin{center}
\begin{tabular}{l}
$m=h_0n+n_1$ \\ $n=h_1n_1+n_2$ \\ $n_1=h_2n_2+n_3$ \\ \ \vdots \\
$n_{s-2}=h_{s-1}n_{s-1}+1$ \\ $n_{s-1}=h_s1$.
\end{tabular}
\end{center}

We denote by $\frac{m}{n}=[h_0,\ldots,h_s]$ the partial fraction decomposition of $\frac{m}{n}$, adjusted in such a way that $s$ becomes  even, say $s=2t$ (for example, $[a_0,a_1]=[a_0,a_1-1,1]$). 
Put $\frac{q_i}{p_i}=[h_0,\ldots,h_{i}]$ in such a way that $q_i$ and $p_i$ are coprime. So, one has the following theorem:

\begin{teor}[Casas-Alvero \cite{C1}] \label{Casas} If $f$ is a general member of $K(n,m)$ where ${\rm \gcd}(n,m)=1$, then $P(f)$ has one Merle package with branches $\gamma_{i,j}$, $i=1,\ldots,t$, $j=1,\ldots,h_{2i}$, having multiplicity $\I(f,X)=p_{2i-1}$  and 
$\I(f,Y)=q_{2i-1}$ and such that $$\I(\gamma_{i,j},\gamma_{i',j'})  =\min(p_{2i-1}q_{2i'-1},p_{2i'-1}q_{2i-1}).$$
\end{teor}

\begin{remark} Notice that the branches of $P(f)$ for a general $f\in K(n,m)$ are all smooth if and only if $p_{2i-1}=1$, for all $i$. But, since the $p_i$ form an increasing sequence, this only may happen when $2t-1=1$, that is, $t=1$. 

If $\frac{m}{n}=[h_0,h_1,h_2]$, then we have $m=h_0n+n_1$; \ $n=h_1n_1+1$; \ $n_1=h_2\cdot 1$. The condition that $\frac{q_1}{p_1}=[h_0,h_1]$ is an integer is equivalent to $h_1=1$ and $h_2=n-1$. Hence the fact that $P(f)$ has only smooth branches is equivalent to $m=(h_0+1) n-1$. 

In the case where $\frac{m}{n}=[h_0,h_1-1,1]$, so $\frac{q_1}{p_1}=[h_0,h_1-1]$. Now, the condition that $\frac{q_1}{p_1}$ is an integer is equivalent to $h_1=2$ and this in turn is equivalent to $n=2$. Hence, the fact $P(f)$ has only smooth branches is equivalent to $m=h_0\cdot 2+1=(h_0+1)2-1$.

In conclusion, one has that $P(f)$, where $f$ corresponds to a general member of $K(n,m)$, has only smooth branches, if and only if $m=\lambda n-1$, where $\lambda$ is some natural number greater than $1$.
\end{remark}

\subsection{The infinitely near points}
Let $S_0 \subset \C^2$ be an open set containing the origin $0=(0,0)$. Let $\pi\colon S_1 \to S_0$ be the blow-up
of $S_0$ centered at $0$ and denote by $E_0=\pi^{-1}(0)$ the exceptional divisor of $\pi$. We denote by $\mathcal{N}_0$ the set of infinitely near points to $0$, which can be viewed as the disjoint union of $0$ and all exceptional
divisors obtained by successive blowing-ups above $0$. The set of points on 
the exceptional divisor of the $i$-th blow-up centered at a point $P\in S_{i-1}$ are called the first
infinitesimal neighborhood of $P$ and the $i$-th infinitesimal neighborhood of $0$.
The set $\mathcal{N}_0$ is naturally endowed with
an order relation defined by $P < Q$ if and only if $Q \in
\mathcal{N}_P$.

Given $f\in \C\{x,y\}$ that defines a curve $C$ and given $P$ in the first infinitesimal neighborhood of $0$, we
denote by $C_P$ the germ of curve at $P$ defined via the strict transform $\widetilde{f}_P$ of $f$, which might be
viewed as the germ at $P$ of the closure of $\pi^{-1}(C\setminus\{0\})$. By induction we
may obtain the strict transform of $C$ at any point of $\mathcal{N}_0$.

The multiplicity of $C_P$ at $P\in \mathcal{N}_0$ is $m_P(f) = m_P(\widetilde{f}_P)$.
We say that $P$ lies on $C$, or belongs to it, if and only if $m_P(f) >
0$, and denote by $\mathcal{N}_0(f)$ the set of all such points. A
point $P \in \mathcal{N}_0(f)$ is {\it simple} (resp. {\it
multiple}) if and only if $m_P(f) = 1$ (resp. $m_P(f) > 1$). Given
two germs of curves $C_f$ and $C_g$, their intersection multiplicity at $0$
can be computed by means of Noether's formula as follows:
\begin{equation}\label{Noether}
\I(f,g)=\sum_{P\in \mathcal{N}_0(f)\inter\mathcal{N}_0(g)}
m_P(f)m_P(g).
\end{equation}

Given $P, Q\in \mathcal N_0$ such that $P<Q$, we say that $Q$ is {\it
proximate} to $P$ (writing $Q\to P$) if and only if $Q$
lies on the exceptional divisor $E_P$ or in the strict
transform of $E_P$. A point $P$ is said to be {\it free} (resp. {\it
satellite}) if it is proximate to exactly one point (resp. two
points), and these are the only possibilities. Notice that
$Q\to P$ implies $Q > P$, but not conversely. 

An important formula due to Noether is the following:
$$m_P(f)=\sum_{Q \to P}m_Q(f).$$

A point $P\in\mathcal{N}_0(f)$ is {\it singular} if it is either multiple,
or satellite, or precedes a satellite point on $C_f$, and it is
{\it non-singular}, or {\it regular}, otherwise. Equivalently, $P$ is
non-singular if and only if it is free and there is no satellite
point $Q > P$.

Let $C_f=\bigcup_{i=1}^{s}C_{f_i}$ be a reducible plane curve we
denote by $P_i\in\mathcal{N}_0(f)$ the first regular point on
$C_{f_i}$. We denote by $$S(f)=\{Q\in\mathcal{N}_0(f);\ Q=P_i\
\mbox{or}\ Q\ \mbox{is singular}\}.$$

It may be shown that 
two curves $C_f$ and $C_g$ are equisingular if and only
if there exists a bijection $\phi$$:S(f)\to S(g)$ such that both $\phi$, $\phi^{-1}$ preserve the natural ordering and the proximity relations among their infinitely  near points.

\begin{defi}
A cluster $\mathcal K$ is a finite
subset $K \subset \mathcal{N}_0$ such that if $P \in K$, then any
other point $Q < P$ also belongs to $K$, together with a valuation $ v_{\mathcal K}\colon K
\longrightarrow \Z$. The set $K$ is called the support of $\mathcal K$ and the number $v_{\mathcal K}(P)$ is the {\it virtual
multiplicity} of $P$ in $\mathcal K$.
\end{defi}

We follow Casas, representing a cluster by means of an Enriques
diagram, which is a  tree whose vertices are identified with the
points in $K$ (the root corresponds to the origin $0$) and there is
an edge between $P$ and $Q$ if and only if $P$ lies on the first
neighborhood of $Q$ or vice-versa. Moreover, the edges are drawn
according to the following rules:
\begin{enumerate}[{i)}]
\item If $Q$ is free and proximate to $P$, the edge joining $P$ and $Q$ is curved and if $P\neq 0$, it is tangent to the edge
ending at $P$.
\item If $P$ and $Q$ ($Q$ in the first neighborhood or $P$) have been represented, the other points proximate
to $P$ in successive neighborhoods of $Q$ are represented on a
straight half-line starting at $Q$ and orthogonal to the edge ending
at $Q$.
\end{enumerate}

\begin{defi}
We will say that a curve $C_f$ goes sharply through the cluster $\mathcal K$ if  $C_f$ goes through $K$ with effective
multiplicities equal to the virtual ones and has no singular points
outside of $K$.
\end{defi}

\subsection{Enriques' Theorem}

In what follows we will describe the {\em cluster of singularities} of a plane branch 
$C_f$, that is, the cluster ${\mathcal K}(f)=(S(f),v_{\mathcal K(f)})$, where $v_{\mathcal K(f)}(P)=m_P(f)$.\smallskip

Suppose that $C_f$ has multiplicity $n$ and characteristic exponents
$m_1,\ldots,m_r$, then $C_f$ is analytically equivalent to a curve that admits a Puiseux
parametrization of the form $x=t^n$, $y=\sum_{i\geq m_1}a_it^i$
such that $a_{m_k}\neq 0$ for $k=1, \ldots ,r$ and
$a_{m_1}=1$.

Denoting $m_{0}=0$,  
$n_0^k=e_{k-1}=\gcd(n,m_1,\ldots,m_{k-1})$, $n_0^{k+1}=n_{s(k)}^k=e_k$, 
we consider the euclidean expansions
\begin{center}
\begin{tabular}{l}\label{L}
$m_k-m_{k-1}=h_0^kn_0^k+n_1^k$ \\
$n_0^k=h_1^kn_1^k+n_2^k$ \\ $n_1^k=h_2^kn_2^k+n_3^k$ \\ \ \vdots \\
$n^k_{s(k)-2}=h^k_{s(k)-1}n^k_{s(k)-1}+n^k_{s(k)}$ \\
$n^k_{s(k)-1}=h^k_{s(k)}n^k_{s(k)}.$
\end{tabular}
\end{center}

When $k=1$, we omit the index $k$ in $n^k_j$, $h^k_j$ and $s(k)$.

The cluster of $C$ is composed by $r$ blocks, which we describe below.

The first block is composed as follows:
\begin{description}
\item[\ ] It starts with the point $P_{0,1}=O$, followed by points $P_{0,i}\in\mathcal{N}_0(f)$, $i=2,\ldots ,h_0$, each one
in the first neighborhood of the preceeding one, all free with value $n$.
\item[\ ] It continues with the point $P_{1,1}$, free in the fist neighborhood of
$P_{0,h_0}$, followed by points $P_{1,i}$, $i=2,\ldots ,h_1$, not free and each in the
first neighborhood of the preceeding one, with value $n_1$.

\item[\ ] For $2\leq j\leq s$, the point $P_{j,1}$ is proximate to $P_{j-2,h_{j-2}}$ and for $i=2,\ldots ,h_j$ we have $P_{j,i}$
proximate to $P_{j-1,h_{j-1}}$ in the first neighborhood of
$P_{j,i-1}$ with value $n_j$.
\end{description}

For $1<k\leq r$, we put $P^{k-1}_{s(k-1),h^{k-1}_{s(k-1)}}=P^k_{0,0}$. The points of the cluster in the $k$-th block after $P^k_{0,0}$ are given by:

\begin{description}
\item[\ ] $h^k_0$ free points $P^k_{0,1},\ldots ,P^k_{0,h^k_0}$ with
value $n^k_{0}$;

\item[\ ] $h^k_1$ points $P^k_{1,1},\ldots ,P^k_{1,h^k_1}$ with
value $n^k_{1}$ proximate to $P^k_{0,h^k_0}$.

\item[\ ]  For $2\leq j\leq s(k)$, we have $h^k_{j}$ points $P^k_{j,1},\ldots ,
P^k_{j,h^k_{j}}$, where the first one is proximate to $P^k_{j-2,h^k_{j-2}}$
and for $i=2,\ldots ,h^k_j$, the point $P^k_{j,i}$ is proximate to
$P^k_{j-1,h^k_{j-1}}$ and all of them have value $n^k_j$.
\end{description}

This yields to the following Enriques diagrams:

{\tiny
\begin{center}
\setlength{\unitlength}{1cm}
\begin{picture}(15,3)
\put(0,0){\qbezier(0,0),(0.2,1),(1,2)}\put(0,0){\circle*{0.1}$P_{0,1}$}\put(0.2,0.7){\circle*{0.1}$P_{0,2}$}\put(0.8,1.7){\circle*{0.1}}
\put(0,1.7){$P_{0,h_0}$}
\put(1,2){\circle*{0.1}$P_{1,1}$}\put(1,2){\line(1,-1){1}}\put(1.2,1.8){\circle*{0.1}}\put(1.7,1.3){\circle*{0.1}}
\put(1,1.1){$P_{1,h_1}$}
\put(2,1){\circle*{0.1}}\put(2,0.8){$P_{2,1}$}\put(2,1){\line(1,1){1}}\put(2.2,1.2){\circle*{0.1}}
\put(2.8,1.8){\circle*{0.1}}\put(2.8,1.6){$P_{2,h_2}$}
\put(3,2){\circle*{0.1}$P_{3,1}$}\put(3.8,1.2){$\ddots$}
\put(5,0.6){\line(1,1){1}}\put(5,0.6){\circle*{0.1}
    $ P^{i-1}_{s(i-1),1}$ }\put(5.8,1.4){\circle*{0.1}}
\put(6,1.6){\circle*{0.1}} \put(6.1,1.5){$P^{i-1}_{s(i-1),h^{i-1}_{s(i-1)}}$}
\put(6,1.6){\qbezier(0,0),(0.2,0.5),(1.8,1)}\put(6.3,1.9){\circle*{0.1}}\put(5.6,1.9){$P^i_{0,1}$}\put(7,2.3){\circle*{0.1}}\put(6.2,2.6){$P^i_{0,h^i_0}$}
\put(7.8,2.6){\circle*{0.1}}\put(7.8,2.8){$P^i_{1,1}$}

\put(7.8,2.6){\line(1,-1){1}}\put(8,2.4){\circle*{0.1}$P^i_{1,2}$}\put(8.6,1.8){\circle*{0.1}}\put(7.7,1.9){$P^i_{1,h^2_1}$}
  \put(8.8,1.6){\circle*{0.1}}\put(8.5,1.1){$P^i_{2,1}$}
\put(8.8,1.6){\line(1,1){1}}\put(9.3,2.1){\circle*{0.1}}\put(9.8,2.6){\circle*{0.1}}
\put(9.6,2.8){$P^i_{3,1}$}\put(9.8,2.3){$\ddots$}\put(10.1,2.1){$\ddots$}
\put(11,1.6){\line(1,1){1}}\put(11,1.6){\circle*{0.1}}\put(10.8,1.2){$P^r_{s(r),1}$}
\put(11.6,2.2){\circle*{0.1}}\put(12,2.6){\circle*{0.1}}\put(11.4,2.8){$P^r_{s(r),h^r_{s(r)}}$}
\end{picture}
\end{center}}
\begin{center}	
{\small if \ $m_i-m_{i-1}>e_{i-1}.$}
\end{center}

{\tiny
\begin{center}
    \setlength{\unitlength}{1cm}
    \begin{picture}(15,3)
    \put(0,0){\qbezier(0,0),(0.2,1),(1,2)}\put(0,0){\circle*{0.1}$P_{0,1}$}\put(0.2,0.7){\circle*{0.1}$P_{0,2}$}\put(0.8,1.7){\circle*{0.1}}
    \put(0,1.7){$P_{0,h_0}$}
    \put(1,2){\circle*{0.1}$P_{1,1}$}\put(1,2){\line(1,-1){1}}\put(1.2,1.8){\circle*{0.1}}\put(1.7,1.3){\circle*{0.1}}
    \put(1,1.1){$P_{1,h_1}$}
    \put(2,1){\circle*{0.1}}\put(2,0.8){$P_{2,1}$}\put(2,1){\line(1,1){1}}\put(2.2,1.2){\circle*{0.1}}
    \put(2.8,1.8){\circle*{0.1}}\put(2.8,1.6){$P_{2,h_2}$}
    \put(3,2){\line(1,-1){1}}\put(3,2){\circle*{0.1}$P_{3,1}$}\put(3.8,1.2){\circle*{0.1}}\put(3,1.2){$P_{3,h_3}$}
    \put(4,1){\circle*{0.1}} \put(4,0.6){$P_{4,1}$}
    \put(4,1){\line(1,1){1}}\put(5,2){\circle*{0.1}}\put(5.2,1.5){$\ddots$}
    \put(6,0.6){\line(1,1){1}}\put(6,0.6){\circle*{0.1}
        $ P^{i-1}_{s(i-1),1}$ }\put(6.8,1.4){\circle*{0.1}}
    \put(7,1.6){\circle*{0.1}} \put(7.1,1.5){$P^{i-1}_{s(i-1),h^{i-1}_{s(i-1)}}$}
    \put(7,1.7){\qbezier(0,0),(0.1,0.2),(0.3,0.2)}
    \put(7.3,1,9){\circle*{0.1}}\put(6.8,2){$P^i_{1,1}$}
    \put(7.3,1.9){\line(2,1){1}}\put(8.3,2.4){\circle*{0.1}}
    \put(8.3,2.5){$P^i_{2,1}$}
    \put(8.3,2.4){\line(1,-1){1}}
    \put(9.3,1.4){\circle*{0.1}}\put(9.3,1.5){$P^i_{2,h^i_2}$}
    \put(9.3,1.5){\line(1,1){1}}
    \put(10.3,2.5){\circle*{0.1}}\put(10.3,2.6){$P^i_{3,1}$}
    \put(10.3,2.2){ $\ddots $}
    \put(11,1.5){\circle*{0.1}}\put(11,1.6){$P^r_{s(r),1}$}
    \put(11,1.5){\line(1,1){1}}
    \put(12,2.5){\circle*{0.1}}
    \put(11.3,2.9){$P^r_{s(r),h^r_{s(r)}}$}
    \end{picture}
   \end{center}}
	\begin{center}
	{\small if \ $m_i-m_{i-1}<e_{i-1}$.}
	\end{center}

\section{Description of the packages in Merle's Theorem}

By Casas' theorem \cite{C2}, we have that the cluster $\mathcal K^r$ of the polar of a branch corresponding to a general member of $K(n,m_1,\ldots,m_r)$ has the same support $K^r$ as the cluster $\mathcal K(f)$ of the singularities of $C_f$, that is
 $$K^r=\{P^k_{i,j}; 1\leq k \leq r, 0\leq i\leq s(k), 1\leq
j\leq h^i_{s(i)}\},$$ with valuation:
\begin{equation}\label{valuations}
v_{\mathcal K^r}(P^k_{i,j})= \left\{
 \begin{array}{ll}
m_{P^k_{s(k),h^k_{s(k)}}}(f)-1, &\ \ \mbox{if $(i,j)=(s(k), h^k_{s(k)})$; \ otherwise},\\ 
m_{P^k_{i,j}}(f)-1,&\ \ \mbox{if $i$\ is even},\\ 
m_{P^k_{i,j}}(f),&\ \ \mbox{if\ $i$\ is odd}.
\end{array}
\right.
\end{equation}

To describe explicitly Merle's packages of such a polar, we firstly consider the cluster $\mathcal K^{'}$ given as follows:
\medskip

\noindent {\bf 1.} \  If $m_r-m_{r-1}> e_{r-1}$, then its support is $K^{'}=\{ P^r_{i,j};  0\leq i\leq s(r), 1\leq
j\leq h^r_{s(r)}\}$, with valuation: 

$$
v_{\mathcal K^{'}}(P^r_{i,j})=\left\{ 
\begin{array}{ll}
0,&\ \ \mbox{if $(i,j)=(s(r),h^r_{s(r)})$}; \mbox{otherwise,}\\ 
n^r_i-1, &\ \ \mbox{if $i$ is even},\\ 
n^r_i,&\ \ \mbox{if $i$ is\ odd}.
\end{array}\right.
$$ 

Notice that $\mathcal K^{'}$ represents the cluster of the polar of a general curve $f_r$ in
$K(e_{r-1},m_r-m_{r-1})$ based at $P^r_{1,0}$.\medskip

\noindent {\bf 2.} \ If $m_r-m_{r-1}< e_{r-1}$, then its support is $$K^{'}= \{ P^r_{i,j};  1\leq i\leq s(r), 1\leq
j\leq h^r_{s(r)}\}\cup \{P^{r-1}_{s(r-1),h^{r-1}_{s(r-1)}}\},$$ with same values as above on the first set and 
$$
v_{\mathcal K^{'}}(P^{r-1}_{s(r-1),h^{r-1}_{s(r-1)}})  =e_{r-1}-1.$$

Notice that this represents the cluster of the polar of a general curve $f_r$ in $K(e_{r-1},m_r-m_{r-1}+e_{r-1})$ based at $P^{r-1}_{s(r-1),h^{r-1}_{s(r-1)}}$.
\medskip

Now, by Theorem \ref{Casas}, we have that:
$$P(f_r) =
\prod_{i=1}^{[\frac{s(r)+1}{2}]}\prod_{j=1}^{h^r_{2i}}\gamma^r_{i,j}$$
where $\gamma^r_{i,j}$ is determined by $m_r-m_{r-1}$ and $e_{r-1}$
according to the following cases:\medskip

\noindent {\bf 1$'$.} If $m_r-m_{r-1} >e_{r-1}$, writing $\frac{m_r-m_{r-1}}{e_{r-1}}=[h^r_0,\ldots ,h^r_{s(r)}]$, one has that
 $\gamma^r_{i,j}\in K(p^r_{2i-1},q^r_{2i-1})$,
where $\frac{q^r_{2i-1}}{p^r_{2i-1}}=[h^r_0,\ldots ,h^r_{2i-1}]$ and $\gcd(p^r_{2i-1},q^r_{2i-1})=1$.\medskip

\noindent {\bf 2$'$.}  If $m_r-m_{r-1} < e_{r-1}$, writing $\frac{m_r-m_{r-1}+e_{r-1}}{e_{r-1}}=[1,h^r_1,\ldots,
h^r_{s(r)}]$, one has that $\gamma^r_{i,j}\in
K(p^r_{2i-1},p^r_{2i-1}+q^r_{2i-1})$, where $\frac{q^r_{2i-1}}{p^r_{2i-1}}=[0,h^r_1,\ldots ,h^r_{2i-1}]$.
\medskip 

Now, by blowing down the branches $\gamma^r_{i,j}$ to the point $P_{0,1}$, with respect to the cluster of singularities of any element in $K(
\widetilde{n},\widetilde{m}_1,\ldots,\widetilde{m}_{r-1})$, where $\widetilde{n}=n/e_{r-1}$ and $\widetilde{m_i}=m_i/e_{r-1}$, $i=1,\ldots,r-1$, we get branches 
$\xi^r_{i,j}$ that pass through the points $K^{r-1}\cup K_i^{'}$, where $K^{r-1}=\{P^k_{i,j}; 1\leq k \leq r-1, 0\leq i\leq s(k), 1\leq
j\leq h^i_{s(i)}\}$, and $K_i^{'}=\{P_{0,1}^r,\ldots,P^r_{2i-1,h_{2i-1}^r}\}$,
with multiplicities at the points of $K^{r-1}$ given by 
$$
m_{P^k_{i,j}}(\xi^r_{i,j})=
\frac{m_{P^k_{i,j}}(f)}{e_{r-1}} p^r_{2i-1}, \ \ k=1,\ldots,r-1;
$$
and the multiplicities at the points of $K^{'}$  given according the following cases:\medskip

\noindent {\bf 1$''$.} \ For $m_r-m_{r-1}>e_{r-1}$ we have that the multiplicities of the
$\xi^r_{i,j}$ at the points $P^r_{0,1},\cdots,
P^r_{2i-1,h^r_{2i-1}}$ are determined by
$\frac{q^r_{2i-1}}{p^r_{2i-1}}$, then by definition of $\xi^r_{i,j}$
we have that the strict transform of the curve $\xi_{r}$ in the
point $P^r_{0,1}$ goes sharply through the cluster
$\mathcal K^{'}$, since the strict transform of $\xi^r_{i,j}$ at the point $P^r_{0,1}$ coincides with $\gamma^r_{i,j}$.\smallskip

\noindent {\bf 2$''$.} \ For $m_r-m_{r-1}<e_{r-1}$ we have that the multiplicities of the $\xi^r_{i,j}$ at the points
$P^{r-1}_{s(r-1),h^{r-1}_{s(r-1)}}, P^r_{0,1},\ldots,
P^r_{2i-1,h^r_{2i-1}}$ are determined by
$1+\frac{q^r_{2i-1}}{p^r_{2i-1}}$ and, from the definition of
$\xi^r_{i,j}$, the strict transform of curve $\xi_{r}$ at the
point $P^{r-1}_{s(r-1),h^{r-1}_{s(r-1)}}$ goes sharply through the
cluster $\mathcal K^{'}$, for the same reason as above.\smallskip

From the above analysis, one sees that 
$$\xi_r=\prod_{i=1}^{[\frac{s(r)+1}{2}]}\prod_{j=1}^{h^r_{2i}}\xi^r_{i,j},$$
with
$$
\xi^r_{i,j}\in K(p^r_{2i-1}\widetilde{n},p^r_{2i-1}\widetilde{m}_1,\ldots,p^r_{2i-1}\widetilde{m}_{r-1},p^r_{2i-1}\widetilde{m}_{r-1}+q^r_{2i-1}).
$$


In order to describe the decomposition of the polar of $f$ we consider the cluster
$\overline{\mathcal K}$ whose support is the same as that of $\mathcal K(f)$ (or of $\mathcal K^r$), with valuation
$v_{\overline{\mathcal K}}(P^k_{i,j})=v_{\mathcal K^r}(P^k_{i,j})-m_{P^k_{i,j}}(\xi_r)$.

In particular, we have that $v_{\overline{\mathcal K}}(P^r_{i,j})=0$ and if $m_r-m_{r-1}<e_{r-1}$, then
$$v_{\overline{\mathcal K}}(P^{r-1}_{r-1,s(r-1)})=
(e_{r-1}-1)-(e_{r-1}-1)=0.$$
By a computation, using the proximity relations, one obtains
$$v_{\overline{\mathcal K}}(P^{r-1}_{i,j})=\left \{
\begin{array}{ll}
(\widetilde{n}^{r-1}_ie_{r-1}-1)-(e_{r-1}-1)\widetilde{n}^{r-1}_i=
\widetilde{n}_i^{r-1}-1, & \mbox{if $i$ is even,}\\ \\
\widetilde{n}^{r-1}_i e_{r-1}-( e_{r-1}-1)\widetilde{n}^{r-1}_i = \widetilde{n}_i^{r-1}, & \mbox{if $i$ is odd}.
\end{array}
 \right .$$
where $\widetilde{n}^{r-1}_i=\frac{n^{r-1}_i}{e_{r-1}}$.\medskip

Using Noether's formulas and by a similar argument, it is possible
to show that $v_{\overline{\mathcal K}}(P^k_{i,j})=\widetilde{n}^k_i-1$, if
$i$ is even, and $v_{\overline{\mathcal K}}(P^k_{i,j})=\widetilde{n}^k_i$, if
$i$ is odd.

Finally, in any situation we have
$v_{\overline{\mathcal K}}(P^k_{s(k),h^k_{s(k)}})=\widetilde{n}^k_{s(k)}-1$.
In this way, the cluster $\overline{\mathcal K}$ represents the cluster of
singularities of the polar curve of a generic branch $g$ in
$K(\widetilde{n},\widetilde{m}_1,\cdots,\widetilde{m}_{r-1})$.
Therefore,
 $$P(f)= P(g)\xi_r.$$

Now, repeating the same procedure to $P(g)$, and so on, we obtain
$$P(f) = \xi_1\cdots \xi_{r-1} \xi_r,$$
where $\xi_1 = P(f_1)$ and $f_1$ is a general member of $
K(\frac{n}{e_1},\frac{m_1}{e_1})$, which is explicitly described in Theorem \ref{Casas}. On the other hand, 
\begin{equation} \label{eq} \xi_{k+1}=\prod_{i=1}^{[\frac{s(k+1)+1}{2}]}\prod_{j=1}^{h_{2i}}\xi^{k+1}_{i,j}, \ \ k=1,\ldots,r-1,\end{equation} where, if we write $\frac{m_{k+1}-m_{k}}{e_{k}}=[h^{k+1}_0,\ldots,h^{k+1}_{s(k+1)}]$ and define
$$\frac{q^{k+1}_{2i-1}}{p^{k+1}_{2i-1}}=[h^{k+1}_0,h^{k+1}_1,\ldots,
h^{k+1}_{2i-1}], \ \text{with} \ \gcd(p^{k+1}_{2i-1},q^{k+1}_{2i-1})=1,$$ we have 
$$\xi^{k+1}_{i,j}\in K(p^{k+1}_{2i-1} \frac{n}{e_k},p^{k+1}_{2i-1}
\frac{m_1}{e_k},\ldots, p^{k+1}_{2i-1} \frac{m_k}{e_k},p^{k+1}_{2i-1}
\frac{m_k}{e_k}+q^{k+1}_{2i-1} ).$$

Summarizing, we have proved part of the following result.

\begin{teor} \label{main} If $f$ is a general branch in $K(n,m_1,\ldots,m_r)$, then the Merle decomposition of $P(f)$ is given by
$$P(f) = \xi_1\xi_2\cdots \xi_r,$$
where $\xi_1 = P(f_1)$ with $f_1$ a general member of 
$K(\frac{n}{e_1},\frac{m_1}{e_1})$ and 
$\xi_{k+1}$ is as in (\ref{eq}).

The intersection multiplicities of these branches are given by:\smallskip

\noindent  $\I(\xi^{k+1}_{i,j},\xi^{k+1}_{u,v})= p^{k+1}_{2i-1}p^{k+1}_{2u-1} \left(\frac{n}{e_k^2}+\sum_{w=1}^{k-1}\frac{e_w}{e_{k}^2}(m_{w+1}-m_w)\right)+ q^{k+1}_{2u-1}p^{k+1}_{2i-1},$ for $i\leq u$.\smallskip

\noindent $\I(\xi^{l+1}_{i,j},\xi^{k+1}_{u,v})=
        \frac{p^{l+1}_{2i-1}p^{k+1}_{2u-1}}{e_{l}e_{k}}\big(\sum_{w=1}^{l}m_w(e_{w-1}-e_w)+m_{l +1}e_{l}\big)$, \  for $k>l\geq 0$, with the convention that $\sum_{w=s}^tA_w=0$, if $t<s$.
						
						\end{teor}
\prova: It remains only to compute the intersection multiplicities.

By  Noether's formula, we know that the
intersection multiplicity of two branches is the sum of the
 products of the multiplicities in common points.\smallskip

\noindent Case 1.\ The branches belong to the same package. 

Suppose that $1\leq i\leq u\leq \left[ \frac{s(k+1)+1}{2}\right]$, $1\leq j \leq h_{2i}$ and $1\leq v \leq h_{2u}$ and let

$$\begin{array}{l}
\xi^{k+1}_{i,j}\in K(p^{k+1}_{2i-1}\frac{n}{e_{k}},p^{k+1}_{2i-1}\frac{m_1}{e_{k}},\ldots,p^{k+1}_{2i-1}\frac{m_i}{e_{k}},p^{k+1}_{2i-1}\frac{m_i}{e_{k}}+q^{k+1}_{2i-1}),  \ \hbox{and}\\ \\
\xi^{k+1}_{u,v}\in K(p^{k+1}_{2u-1}\frac{n}{e_{k}},p^{k+1}_{2u-1}\frac{m_1}{e_{k}},\ldots,p^{k+1}_{2u-1}\frac{m_i}{e_{k}},p^{k+1}_{2u-1}\frac{m_i}{e_{k}}+q^{k+1}_{2u-1}  ).
\end{array}
$$

As $i \leq u $, we have that the last common point of the two above branches is $P^{k+1}_{2i-1,h^{k+1}_{2i-1}}$. 
Using the clusters of both branches, we obtain that the sum of products of the multiplicities until the point $P^{k}_{s(k),h^{k}_{s(k)}} $ is $$\left( \frac{e_1}{e_{k}^2}m_1+\sum_{j=1}^{k-1}\frac{e_j}{e_{k}^2}(m_{j+1}-m_j)\right)p^{k+1}_{2i-1}p^{k+1}_{2u-1}.$$ 

On the other hand, since the branches at the point $P^{k+1}_{0,1}$ are the branches of the polar of a genus one curve, using Theorem \ref{Casas}, one gets
$$\I_{P^{k+1}_{0,1}}(\xi^{k+1}_{i,j},\xi^{k+1}_{u,v})=q^{k+1}_{2u-1}p^{k+1}_{2i-1}.$$ 

Summing up and using Noether's formula, one gets that
$$
\I(\xi^{k+1}_{i,j},\xi^{k+1}_{u,v})= \left( \frac{e_j}{e_k^2}+ \sum_{w=1}^{k-1}\frac{e_j}{e_{k}^2}(m_{w+1}-m_w)\right)p^{k+1}_{2i-1}p^{k+1}_{2u-1} + q^{k+1}_{2u-1}p^{k+1}_{2i-1}.$$

\noindent Case 2. \ The branches are in distinct packages.

Consider $\xi^{l+1}_{i,j}$ and $\xi^{k+1}_{u,v}$ where $0 \leq
l<k$, $1 \leq i \leq [\frac{s(l+1)+1}{2}]$, $1 \leq u \leq
[\frac{s(k+1)+1}{2}]$, $ 1\leq j \leq h^{l+1}_{2i}$ and $1\leq v
\leq h^{k+1}_{2u}.$

We have that the sum of products of the multiplicities until the point
$P^{l}_{s(l),h^{l}_{s(l)}}$ is
$$\frac{p^{l+1}_{2i-1}p^{k+1}_{2u-1}}{e_{l}e_{k}}\left( nm_1+ \sum_{w=1}^{l-1}e_w(m_{w+1}-m_w)\right),$$
while the sum of products of the multiplicities at the remaining points is
$$\frac{p^{l+1}_{2i-1}p^{k+1}_{2u-1}}{e_{l}e_{k}}(m_{l+1}-m_{l})e_{l}.$$

Therefore, if $e_0=n$ and $m_0=0$, then 
$$\begin{array}{lcl}
\I(\xi^{l+1}_{i,j},\xi^{k+1}_{u,v}) &=& 
\frac{p^{l+1}_{2i-1}p^{k+1}_{2u-1}}{e_{l}e_{k}}\left(nm_1+\sum_{w=1}^{l}e_w(m_{w+1}-m_w) \right)\\
&=&  \frac{p^{l+1}_{2i-1}p^{k+1}_{2u-1}}{e_{l}e_{k}}\big(\sum_{w=1}^{l}m_w(e_{w-1}-e_w)+m_{l
+1}e_{l}\big). 
\end{array}
$$

By construction and by an
analogous computation, we may show that
$$\frac{\I(\xi^{l+1}_{i,j},f)}{m(\xi^{l+1}_{i,j})}=
\frac{1}{n}\big(\sum_{w=1}^{l}m_w(e_{w-1}-e_w)+m_{l
+1}e_{l}\big).$$

In this way, we see that $\xi_u$
is precisely the $u$-th package in Merle's Theorem. \fin

From the above theorem we get immediately the following result:

\begin{cor} The number of branches of the $j$-th package $\xi_j$ in Merle's decomposition of the polar of a general member of $K(n,m_1,\ldots,m_r)$ is given by
 $$\sum_{k=1}^{[\frac{s(j)+1}{2}]} h^j_{2k},$$
where the numbers that appear in the formula are obtained from the euclidean divisions described in (\ref{L}).
\end{cor}

\begin{ejm}
    Let $f$ be general member of $K(8,12,14,15)$. The Euclidean divisions in this case are:

    \begin{center}

        \begin{tabular}{|l|l|l|}
            \hline
            $m_1=12\ \text{and}\ n=8$ & $m_2-m_1=2\ \text{and}\ e_1=4$ & $m_3-m_2=1\ \text{and}\ e_2=2$ \\
            \hline
            $12=1(8)+4$ & $2=0(4)+2$&$1=0(2)+1$\\
            $8=2(4)$ & $4=2(2)$& $2=2(1)$\\
            \hline
        \end{tabular}
    \end{center}
    In this way, we have

    \begin{center}
        \setlength{\unitlength}{1cm}

        \begin{picture}(15,3)
        \put(0,0){\qbezier(0,1),(0.2,2),(1,3)}\put(0,1)
        {\circle*{0.1}$8$}

        \put(1,3){\circle*{0.1}$4$}
        \put(1,3){\line(1,-1){1}}
        \put(2,2){\circle*{0.1}}\put(2,1.6){$4$}

        \put(2,2){\qbezier(0,0),(0.4,-0.3),(0.7,0.1)}
        \put(2.75,2.1){\circle*{0.1}}\put(2.75,2.2){$2$}
        \put(2.7,2){\line(1,1){1}}
        \put(3.7,3){\qbezier(0,0),(0.3,0.2),(0.6,0)}
        \put(3.7,3){\circle*{0.1}}\put(3.7,3.1){$2$}
        \put(4.3,3){\line(1,-1){1}}
        \put(4.3,3){\circle*{0.1}$1$}
        \put(5.3,2){\circle*{0.1}$1$}
        \put(1.2,0.5){Enriques' diagram of $f$}

        \put(7,0){\qbezier(0,1),(0.2,2),(1,3)}\put(7,1)
        {\circle*{0.1}$7$}

        \put(8,3){\circle*{0.1}$4$}
        \put(8,3){\line(1,-1){1}}
        \put(9,2){\circle*{0.1}}\put(9,1.6){$3$}

        \put(9,2){\qbezier(0,0),(0.4,-0.3),(0.7,0.1)}
        \put(9.75,2.1){\circle*{0.1}}\put(9.75,2.2){$2$}
        \put(9.7,2){\line(1,1){1}}
        \put(10.7,3){\qbezier(0,0),(0.3,0.2),(0.6,0)}
        \put(10.7,3){\circle*{0.1}}\put(10.7,3.1){$1$}
        \put(11.3,3){\line(1,-1){1}}
        \put(11.3,3){\circle*{0.1}}\put(11.3,3.1){$1$}
        \put(12.3,2){\circle*{0.1}}\put(12.3,2.1){$0$}
        \put(7.75,0.5){Enriques' diagram of $P(f)$}

        \end{picture}

    \end{center}
    Since $h^3_0=0$, $h^3_1=1$ and $h^3_2=1$, according to Theorem \ref{main}, the  third package $\xi_3$ of $P(f)$ has just one branch $\xi^3_{1,1} \in K(4,6,7)$ whose Enriques' diagram is \bigskip

    \setlength{\unitlength}{1cm}

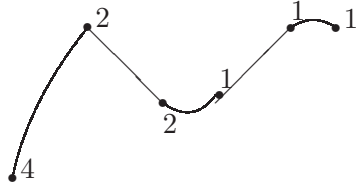
    \begin{center}
        \begin{picture}(7,3)

        \put(0,0){\qbezier(0,1),(0.2,2),(1,3)}\put(0,1)
        {\circle*{0.1}$4$}
                \put(1,3){\circle*{0.1}$2$}
        \put(1,3){\line(1,-1){1}}
        \put(2,2){\circle*{0.1}}\put(2,1.6){$2$}

        \put(2,2){\qbezier(0,0),(0.4,-0.3),(0.7,0.1)}
        \put(2.75,2.1){\circle*{0.1}}\put(2.75,2.2){$1$}
        \put(2.7,2){\line(1,1){1}}
        \put(3.7,3){\qbezier(0,0),(0.3,0.2),(0.6,0)}
        \put(3.7,3){\circle*{0.1}}\put(3.7,3.1){$1$}
                \put(4.3,3){\circle*{0.1}$1$}
        \end{picture}
    \end{center}
\vspace{-1cm}

Now, since $h^2_0=0$, $h^2_1=1$ and $h^2_2=1$, the second package $\xi_2$ of $P(f)$ has just one branch $\xi^2_{1,1}\in K(2,3)$, whose Enriques' diagram is\bigskip

    \setlength{\unitlength}{1cm}
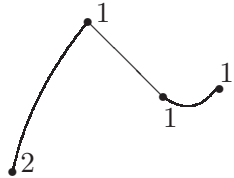
\begin{center}
    \begin{picture}(5,3)

    \put(0,0){\qbezier(0,1),(0.2,2),(1,3)}\put(0,1)
    {\circle*{0.1}$2$}
        \put(1,3){\circle*{0.1}$1$}
    \put(1,3){\line(1,-1){1}}
    \put(2,2){\circle*{0.1}}\put(2,1.6){$1$}
        \put(2,2){\qbezier(0,0),(0.4,-0.3),(0.7,0.1)}
    \put(2.75,2.1){\circle*{0.1}}\put(2.75,2.2){$1$}
    \end{picture}
\end{center}
\vspace{-1cm}

Finally, the first package is $\xi_1$, corresponding to the polar of a general member of $K(2,3)$, hence it has one smooth branch $\xi_{1,1}$, whose Enriques' diagram is \bigskip

    \setlength{\unitlength}{1cm}
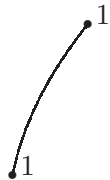
    \begin{center}
        \begin{picture}(5,3)
        \put(0,0){\qbezier(0,1),(0.2,2),(1,3)}\put(0,1)
        {\circle*{0.1}$1$}
        \put(1,3){\circle*{0.1}$1$}
            \end{picture}
    \end{center}
\vspace{-1cm}
		
For the intersection multiplicities of these branches, the theorem gives us		
	$$\I(\xi_{1,1},\xi^2_{1,1})=3,\, \I(\xi_{1,1},\xi^3_{1,1})=6,\,
\I(\xi^2_{1,1},\xi^3_{1,1})=13.$$	
\end{ejm}				
				
From Merle's Theorem it follows that each branch of the $j$-th Merle's package $\xi_j$ of the polar of an irreducible curve  has genus at least $j-1$. On the other hand, from the proof of Theorem \ref{main} one may see that the genus of each component of $\xi_j$ is less or equal than $j$, when the curve is general in its equisingularity class. This generality condition is a sufficient condition to guarantee the bound $j$ from above for the genus of the components of $\xi_j$, as one may see in Remark 2.1 of \cite{HHI2}. 

The problem we address now is to characterize the equisingularity classes given by $K(n,m_1,\ldots,m_r)$ for which the general member has its polar curve composed by branches with genus up to $r-1$.

\begin{cor} Let $f$ be a power series corresponding to a general member of $K(n,m_1,\ldots,m_r)$. The polar of $f$ has branches of genus at most $r-1$, if and only if $m_r=m_{r-1}+\lambda e_{r-1}-1$, for some integer $\lambda\geq 1$.
\end{cor}
\prova: \ From Theorem 3.1, this happens if and only if the $\xi^{r}_{i,j}$ have genus $r-1$. Since $\xi^{r}_{i,j}\in K(p^{r}_{2i-1} \frac{n}{e_{r-1}},p^{r}_{2i-1}\frac{m_1}{e_{r-1}},\ldots, p^{r}_{2i-1} \frac{m_{r-1}}{e_{r-1}},p^{r}_{2i-1}\frac{m_{r-1}}{e_{r-1}}+q^{r}_{2i-1} )$, this, in turn, happens if and only if $p^r_{2i-1}=1$ for all $i=1,\ldots, t(r)$, where $s(r)=2t(r)$. Now, since the  $p^r_{j}$ form an increasing sequence, one must have $t(r)=1$. We have two possibilities:

\noindent 1) \  $m_r-m_{r-1}=h^r_0e_{r-1}+n^r_1$, $e_{r-1}=h^r_1n^r_1+1$ and $n^r_1=h^r_2\cdot 1$. Now, since $\frac{q^r_1}{p^r_1}=[h^r_0,h^r_1]$ is an integer, we must have $h^r_1=1$. Therefore, the condition that $P(f)$ has branches of genus at most $r-1$ is equivalent to 
$$
m_r-m_{r-1}=(h^r_0+1)e_{r-1}-1.
$$

\noindent 2) \ $m_r-m_{r-1}=h^r_0e_{r-1}+1$ and $e_{r-1}=(h^r_1-1) \cdot 1 +1$. Since $\frac{q^r_1}{p^r_1}=[h^r_0,h^r_1-1 ]$ is an integer, then $h^r_1=2$. Which gives  $e_{r-1}=2$. Therefore, the condition that $P(f)$ has branches of genus at most $r-1$ is equivalent to
$$m_r-m_{r-1}=(h^r_0+1)e_{r-1}-1.$$

Concluding in this way our proof.
\fin

\end{document}